\documentclass[onecolumn,draft]{IEEEtran}

\usepackage{enumitem,kantlipsum}
\usepackage{graphicx}
\usepackage{epstopdf}
\usepackage{amssymb}
\usepackage{amsmath}
\usepackage{subfigure}
\usepackage{epsfig}
\usepackage{color}
\usepackage{enumitem}
\usepackage{float}
\usepackage{soul}
\usepackage{wrapfig}
\usepackage{arydshln}
\usepackage{tablefootnote}

\usepackage{amsthm}

\def\0{{\bf 0}}
\def\1{{\bf 1}}

\def\beq{\begin{equation*}}
    \def\eeq{\end{equation*}}
\def\bql{\begin{equation}}
    \def\eql{\end{equation}}
\def\bqn{\begin{eqnarray*}}
    \def\eqn{\end{eqnarray*}}
\def\bnl{\begin{eqnarray}}
    \def\enl{\end{eqnarray}}
\def\bma{\begin{bmatrix}}
    \def\ema{\end{bmatrix}}
\def\bmx{\begin{matrix}}
    \def\emx{\end{matrix}}
\def\ben{\begin{enumerate}}
    \def\een{\end{enumerate}}
\def\bit{\begin{itemize}}
    \def\eit{\end{itemize}}
\def\bei{\begin{itemize}}
    \def\eei{\end{itemize}}
\def\bet{\begin{tabular}}
    \def\eet{\end{tabular}}

\newcommand{\ba}{\mathbf{a}}

\newcommand{\be}{\mathbf{e}}

\newcommand{\mnorm}[1]{{\left\vert\kern-0.25ex\left\vert\kern-0.25ex\left\vert #1   \right\vert\kern-0.25ex\right\vert\kern-0.25ex\right\vert}}

\usepackage[dvipsnames]{xcolor}

\def\1{{\bf1}}

\def\bit{\begin{itemize}}
\def\eit{\end{itemize}}

\def\be{\begin{equation}}
\def\ee{\end{equation}}
\def\ba{\begin{eqnarray}}
\def\ea{\end{eqnarray}}

\def\bes{\begin{equation*}}
\def\ees{\end{equation*}}
\def\bas{\begin{eqnarray*}}
\def\eas{\end{eqnarray*}}

\usepackage{url}
\usepackage{verbatim}

\newtheorem{Remark 1}{Remark}
\newtheorem{Remark 2}[Remark 1]{Remark}
\newtheorem{Remark 3}[Remark 1]{Remark}
\newtheorem{Remark 4}[Remark 1]{Remark}
\newtheorem{Remark 5}[Remark 1]{Remark}
\newtheorem{Remark 6}[Remark 1]{Remark}
\newtheorem{Remark 7}[Remark 1]{Remark}

\newtheorem{Lemma 1}{Lemma}
\newtheorem{Lemma 2}[Lemma 1]{Lemma}
\newtheorem{Lemma 3}[Lemma 1]{Lemma}
\newtheorem{Lemma 4}[Lemma 1]{Lemma}
\newtheorem{Lemma 5}[Lemma 1]{Lemma}
\newtheorem{Lemma 6}[Lemma 1]{Lemma}
\newtheorem{Lemma 7}[Lemma 1]{Lemma}

\newtheorem{Assumption 1}{Assumption}
\newtheorem{Assumption 2}[Assumption 1]{Assumption}
\newtheorem{Assumption 3}[Assumption 1]{Assumption}
\newtheorem{Assumption 4}[Assumption 1]{Assumption}
\newtheorem{Definition 1}{Definition}
\newtheorem{Theorem 1}{Theorem}
\newtheorem{Theorem 2}[Theorem 1]{Theorem}
\newtheorem{Theorem 3}[Theorem 1]{Theorem}
\newtheorem{Theorem 4}[Theorem 1]{Theorem}
\newtheorem{Theorem 5}[Theorem 1]{Theorem}
\newtheorem{Theorem 6}[Theorem 1]{Theorem}
\newtheorem{Theorem 7}[Theorem 1]{Theorem}
\newtheorem{Theorem 8}[Theorem 1]{Theorem}
\newtheorem{Theorem 9}[Theorem 1]{Theorem}
\newtheorem{Theorem 10}[Theorem 1]{Theorem}

\newtheorem{Proposition 1}{Proposition}
\newtheorem{Proposition 2}[Proposition 1]{Proposition}

\title{\LARGE \bf
The spectral radius of a square matrix can be approximated by its ``weighted" spectral norm}

\author{Yongqiang Wang
\thanks{ The work   was supported in part by the National Science Foundation under Grants  ECCS-1912702, CCF-2106293, CCF-2215088, and CNS-2219487.}
\thanks{Yongqiang Wang is with the Department of Electrical and Computer Engineering, Clemson University, Clemson, SC 29634, USA
{\tt\small{yongqiw}@clemson.edu}
}%

  }

\begin{document}

\maketitle
\thispagestyle{empty}
\pagestyle{empty}

\begin{abstract}
In distributed optimization or Nash-equilibrium seeking over directed graphs, it is crucial to find a matrix norm under which the disagreement of individual agents' states contracts. In existing results, the matrix norm  is usually defined by approximating the spectral radius of the matrix, which is possible when the matrix is real and has zero row-sums.    In this brief note, we show that this technique can be applied to general square complex matrices. More specifically, we prove that the spectral radius of any complex square matrix can be approximated by its ``weighted" spectral norm to an arbitrary degree of accuracy.
\end{abstract}

\section{Introduction}

In distributed optimization over directed graphs (particularly in the presence of noise), one commonly used proof technique is to prove that the disagreement of individual agents' states satisfies the contraction property under a certain non-Euclidian matrix norm \cite{pu2020push,wang2022gradient}. This norm is dependent on the interaction pattern of the agents and is usually defined based on the following properties:
\begin{enumerate}
   \item  the spectral radius of the matrix is less than one;
  \item   when the matrix is real and has zero row-sums, one can find a matrix norm to approximate the spectral radius.
\end{enumerate}
 Please see \cite{xin2019distributed} for an instantiation of this norm in distributed optimization over directed graphs. Finding such a  norm is also important in the convergence analysis of distributed Nash-equilibrium seeking algorithms over directed graphs, particularly in the presence of noises (see, e.g., \cite{wang2022ensuring}).

 In this note, we prove that such a norm exists for general square complex matrices. We also prove that this norm is only dependent on the unitary Schur-decomposition matrix (can be selected as the eigenvectors corresponding to the Jordan-normal-form transformation  in the real matrix case) of the matrix, making it possible to define this norm even when the eigenvalues of the interaction matrix are time-varying. In addition, the matrix norm is a ``weighted" spectral norm, implying that it has an associated inner-product.

{\it Notations}: We use $\mnorm{\bullet}$ to denote   matrix norm and  $\|\bullet\|$ for vector norm. For   a complex square matrix $A\in\mathbb{C}^{n\times n}$, we represent its eigenvalues as $\{\lambda_1,\lambda_2,\cdots,\lambda_n\}$ and its spectral radius as $\rho(A)=\max_{1\leq i\leq n}\{|\lambda_i|\}$.  The spectral norm of   $A$ is defined as $\mnorm{A}_2=\sqrt{\rho(AA^\ast)}=\sup_{x\in \mathbb{C}^n,x\neq 0}\frac{\|Ax\|_2}{\|x\|_2}$, where $\|\bullet\|_2$ is the $l_2$ (or Euclidean) vector norm.


\section{Main Results}
\begin{Theorem 1}
  For any   complex square matrix $A\in\mathbb{C}^{n\times n}$ with spectral radius $\rho(A)$, there always exists a matrix norm $\mnorm{A}$ such that the following two statements  hold simultaneously:
\begin{enumerate}
   \item  $\rho(A)  \leq  \mnorm{A} \leq \rho(A)+\epsilon$ holds for any $\epsilon>0$;
  \item the matrix norm $\mnorm{\bullet}$ is induced by a vector norm $\|\bullet\|$ which satisfies $\|x\|=\|\tilde{A}x\|_2$ for any vector $x\in\mathbb{C}^n$, where $\tilde{A}\in\mathbb{C}^{n\times n}$ is an invertible matrix dependent on $A$.
\end{enumerate}
\end{Theorem 1}
 \begin{proof}

According to Lemma 5.6.10 of \cite{horn2012matrix}, there is a unitary Schur-decomposition matrix $U$ and an upper triangular matrix $\Delta$ such that $A=U^\ast\Delta U$ holds (note that for a real matrix $A$, the eigenvector set corresponding to the Jordan-normal-form transformation matrix can   be used here). Setting $D_t\triangleq {\rm diag}(t,t^2,\cdots,t^n)$ with $t$ a positive scalar, we can obtain the following relationship (see  Lemma 5.6.10 of \cite{horn2012matrix}):
\[
\tilde{\Delta}\triangleq D_t\Delta D_t^{-1}=\left[ \begin{array}{cccccc}\lambda_1&t^{-1}d_{1,2}& t^{-2}d_{1,3}&\cdots&t^{-(n-2)}d_{1,n-1}&t^{-(n-1)}d_{1,n}\\
0&\lambda_2&t^{-1}d_{2,3}&\cdots&t^{-(n-3)}d_{2,n-1}&t^{-(n-2)}d_{2,n}\\
0&0&\lambda_3&\cdots&t^{-(n-4)}d_{3,n-1}&t^{-(n-3)}d_{3,n}\\
\vdots&\ddots&\ddots&\ddots&\ddots&\vdots\\
0&0&\cdots &0&\lambda_{n-1}&t^{-1}d_{n-1,n}\\
0&0&0&\cdots&0&\lambda_n
\end{array}\right],
\]
where $d_{i,j}$ denotes the $(i,j)$th entry of $\Delta$ and $\lambda_1,\lambda_2,\cdots,\lambda_n$ are the eigenvalues of $A$.


Define $\Lambda={\rm diag}(\lambda_1,\lambda_2,\cdots,\lambda_n)$. One can obtain
\[
\begin{aligned}
&\tilde{\Delta} -\Lambda =\delta\tilde{\Delta},
\end{aligned}
\]
with
\[
\delta\tilde{\Delta}=\left[ \begin{array}{cccccc}0&t^{-1}d_{1,2}& t^{-2}d_{1,3}&\cdots&t^{-(n-2)}d_{1,n-1}&t^{-(n-1)}d_{1,n}\\
0&0&t^{-1}d_{2,3}&\cdots&t^{-(n-3)}d_{2,n-1}&t^{-(n-2)}d_{2,n}\\
0&0&0&\cdots&t^{-(n-4)}d_{3,n-1}&t^{-(n-3)}d_{3,n}\\
\vdots&\ddots&\ddots&\ddots&\ddots&\vdots\\
0&0&\cdots &0&0&t^{-1}d_{n-1,n}\\
0&0&0&\cdots&0&0
\end{array}\right].
\]


Since   all non-zero entries of the matrix $\delta\tilde{\Delta}$  approach zero at least linearly with an increase in
$t$, for large enough $t$, we can make the spectral norm of $\delta\tilde{\Delta}$ less than $\epsilon$, i.e.,  $\mnorm{\delta\tilde{\Delta}}_2<\epsilon$, which further implies
\[
\begin{aligned}
\mnorm{\tilde{\Delta}}_2-\mnorm{\Lambda}_2&=\mnorm{\tilde{\Delta} -\Lambda +\Lambda }_2-\mnorm{\Lambda }_2 \leq \mnorm{\tilde{\Delta} -\Lambda }_2+\mnorm{\Lambda }_2-\mnorm{\Lambda }_2 <\epsilon.
\end{aligned}
\]
In a similar way, we can also prove
$
\mnorm{\Lambda }_2-\mnorm{\tilde{\Delta} }_2  <\epsilon
$.

Therefore, we have that the distance between   $\mnorm{\Lambda}_2$ and $\mnorm{\tilde{\Delta}}_2$ is less than any $\epsilon>0$. Given that $\Lambda$ is a diagonal matrix, and hence, its spectral norm $\mnorm{\Lambda}_2$ is equal to $\rho(A)$, the spectral radius of $A$, we have that
  the distance between $\mnorm{\tilde{\Delta}}_2$ and $\rho(A)$ is less than $\epsilon$.

Therefore, by defining the matrix norm as
\begin{equation}\label{eq:norm_defintion}
\mnorm{A}\triangleq \mnorm{D_tU AU^\ast D_t^{-1}}_2=\mnorm{(U^\ast D_t^{-1})^{-1} A(U^\ast D_t^{-1})}_2,
\end{equation} we can construct a matrix norm $\mnorm{A}$ such that $\mnorm{A}\leq \rho(A)+\epsilon$ holds for any matrix $A\in\mathbb{C}^{n\times n}$. Given that $\mnorm{A}\geq\rho(A)$ holds for any matrix norm, we arrive at the first statement.

To prove the second statement, we recall the definition of spectral norm (induced $l_2$-norm) of a matrix
\[
\mnorm{A}_2\triangleq \sup_{x\in \mathbb{C}^n,x\neq 0}\frac{\|Ax\|_2}{\|x\|_2},
\]
 where $\|\bullet\|_2$ represents   the $l_2$ (or Euclidean) vector norm.
Hence, (\ref{eq:norm_defintion}) implies
\[
\mnorm{A}\triangleq \sup_{x\in \mathbb{C}^n,x\neq 0}\frac{\| (U^\ast D_t^{-1})^{-1} A(U^\ast D_t^{-1})x\|_2}{\|x\|_2}= \sup_{y\in \mathbb{C}^n,y\neq 0}\frac{\| (U^\ast D_t^{-1})^{-1} Ay\|_2}{\| (U^\ast D_t^{-1})^{-1}y\|_2},
\]
where we have defined $y\triangleq (U^\ast D_t^{-1})x$. Therefore, by defining the norm of a vector $x\in\mathbb{C}^n$ as $\|x\|=\|(U^\ast D_t^{-1})^{-1}x\|_2$, we arrive at the second statement.
\end{proof}

\begin{Remark 1}
  When   $A$ is a real matrix, the $U$ matrix in the proof of Theorem 1 can be selected as the set of eigenvectors that correspond to the Jordan-normal-form transformation.
\end{Remark 1}
 \section*{Acknowledgement}
The author would like to thanks Shi Pu for  helpful discussions.

\bibliographystyle{IEEEtran}

\bibliography{reference1}
\vspace{-1.2cm}

\end{document}